\documentclass[12pt]{article}

\usepackage{amsmath}
\usepackage{amssymb}
\usepackage{amsfonts}
\usepackage{lscape}

\begin{document}

\fontsize{12}{6mm}\selectfont
\setlength{\baselineskip}{2em}

$~$\\[.35in]
\newcommand{\dss}{\displaystyle}
\newcommand{\raro}{\rightarrow}
\newcommand{\be}{\begin{equation}}

\def\sech{\mbox{\rm sech}}
\def\sn{\mbox{\rm sn}}
\def\dn{\mbox{\rm dn}}
\def\coth{\mbox{\rm coth}}
\def\csch{\mbox{\rm csch}}
\def\arctanh{\mbox{\rm arctanh}}
\thispagestyle{empty}

\begin{center}
{\Large\bf An Intrinsic Characterization }  \\    [2mm]
{\Large\bf of Bonnet Surfaces}   \\   [2mm]
{\Large\bf Based on a Closed Differential Ideal}   \\
\end{center}

\vspace{1cm}
\begin{center}
{\bf Paul Bracken}                        \\
{\bf Department of Mathematics,} \\
{\bf University of Texas,} \\
{\bf Edinburg, TX  }  \\
{78540}   \\
{bracken@utpa.edu}
\end{center}

\vspace{3cm}
\begin{abstract}
The structure equations for a surface are introduced and
two required results based on the Codazzi
equations are obtained from them. Important
theorems pertaining to isometric surfaces are stated 
and a theorem due to Bonnet is obtained. A transformation
formula for the connection forms is developed. It is
proved that the angle of deformation must be harmonic. 
It is shown that the differentials of many of the
important variables generate a closed 
differential ideal. This implies that a coordinate
system exists in which many of the variables
satisfy particular ordinary differential equations,
and these results can be used to characterize Bonnet
surfaces.
\end{abstract}

\vspace{2mm}
MSCs:  53A05, 58A10, 53B05
\vspace{2mm}

\newpage
\section{Introduction}
\numberwithin{equation}{section}

Bonnet surfaces in three dimensional Euclidean space have been
of great interest for a number of reasons as a type of surface {\bf [1-2]}
for a long time. 
Bonnet surfaces are of non-constant mean curvature that
admit infinitely many non-trivial and geometrically distinct isometries 
which preserve the mean curvature function.
Non-trivial isometries are ones that do not extend to
isometries of the whole space $E^3$. Considerable interest
has resulted from the fact that the differential equations
that describe the Gauss equations are classified by the type of 
related Painlev\'e equations they correspond to and they are 
integrated in terms of certain hypergeometric transcendents {\bf [3-5]}.
Here the approach first given by S. S. Chern {\bf [6]}
to Bonnet surfaces is considered.  
The developement is accessible with many new proofs given.
The main intention is to end by deriving an intrinsic 
characterization of these surfaces which indicates they are
analytic. Moreover, it is shown that a type of Lax pair 
can be given for these surfaces and integrated. Several
of the more important functions such as the mean curvature
are seen to satisfy nontrivial ordinary differential equations.

Quite a lot is known about these surfaces. With many results
the analysis is local and takes place under the assumptions 
that the surfaces contain no umbilic points and no critical
points of the mean curvature function. The approach here allows
the elimination of many assumptions and it is found the results 
are not too different from the known local ones. The statements and
proofs have been given in great detail in order to 
help illustrate and display the
interconnectedness of the ideas and results involved.

To establish some information about what is known, consider an oriented, connected,
smooth open surface $M$ in $E^3$ with non-constant mean curvature 
function $H$. Moreover, $M$ admits infinitely many non-trivial and
geometrically distinct isometries preserving $H$. Suppose $U$ is the 
set of umbilic points of $M$ and $V$ the set of critical points of
$H$. Many global facts are known with regard to $U, V$ and $H$,
and a few will now be mentioned.
The set $U$ consists of isolated points, even if there exists only
one non-trivial isometry preserving the mean curvature, moreover
$U \subset V$ {\bf [7-8]}. Interestingly, there is no point in $V-U$ at which
all order derivatives of $H$ are zero, and $V$ cannot contain
any curve segment. If the function by which a nontrivial isometry
preserving the mean curvature rotates the principal frame is considered,
as when there are infinitely many isometries,
this function is a global function on $M$ continuously defined {\bf [9-11]}
As first noted by Chern {\bf [6]}, this function is harmonic. The analysis will begin 
by formulating the structure equations for two-dimensional manifolds.

\section{Structure Equations}

Over $M$ there exists a well defined field of
orthonormal frames which is written as ${\bf x}$, $e_1, e_2,e_3$
such that $x \in M$, $e_3$ is the unit normal at ${\bf x}$ and
$e_1, e_2$ are along principal directions {\bf [12]}. The fundamental
equations for $M$ have the form
\be
d {\bf x} = \omega_1 e_1 + \omega_2 \, e_2,
\quad
d e_1 = \omega_{12} e_2 + \omega_{13} e_3,
\quad
d e_2 = - \omega_{12} e_1 + \omega_{23} e_3,
\quad
d e_3 = - \omega_{13} e_1 - \omega_{23} e_2.
\label{eqII1}
\end{equation}
Differentiating each of these equations in turn, results in a
large system of equations for the exterior derivatives of the
$\omega_i$ and the $\omega_{ij}$, as well as a final equation
which relates some of the forms ${\bf [13]}$. This choice of frame
and Cartan's lemma allows for the introduction of the two
principal curvatures which are denoted by $a$ and $c$ at ${\bf x}$
by writing
\be
\omega_{12} = h \omega_1 + k \omega_2, 
\qquad
\omega_{13} =a \omega_1,
\qquad
\omega_{23} = c \omega_2.
\label{eqII2}
\end{equation}
Suppose that $a >c$ in the following. The mean curvature of $M$
is denoted by $H$ and the Gaussian curvature by $K$. They are
related to $a$ and $c$ as follows 
\be
H = \frac{1}{2} (a + c),
\qquad
K = a \cdot c.
\label{eqII3}
\end{equation}
The forms which appear in \eqref{eqII1} satisfy the fundamental
structure equations which are summarized here {\bf [14]},
\be
\begin{array}{ccc}
d \omega_1 = \omega_{12} \wedge \omega_2, &    &  d \omega_2 = \omega_1 \wedge \omega_{12}  \\
           &    &       \\
d \omega_{13} = \omega_{12} \wedge \omega_{23}  & &  d \omega_{23} = \omega_{13} \wedge \omega_{12}, \\
           &    &        \\
   &  d \omega_{12} = a c \, \omega_2 \wedge \omega_1 =- K \, \omega_1 \wedge \omega_2.  &    \\
\end{array}
\label{eqII4}
\end{equation}
The second pair of equations of \eqref{eqII4} 
are referred to as the Codazzi equation and the
last equation is the Gauss equation.

Exterior differentiation of the two Codazzi equations yields
\be
( da - (a-c) h \omega_2 ) \wedge \omega_1 =0,
\qquad
( d c - (a-c) k \omega_1 ) \wedge \omega_2 =0.
\label{eqII5}
\end{equation}
Cartan's lemma can be applied to the equations in \eqref{eqII5}.
Thus there exist two functions $u$ and $v$ such that
\be
\frac{1}{a-c} \, da - h \omega_2 = (u-k) \omega_1,
\qquad
\frac{1}{a-c} \, dc - k \omega_1 = (v-h) \omega_2.
\label{eqII6}
\end{equation}
Subtracting the pair of equations in \eqref{eqII6} gives an
expression for $d \log (a-c)$
\be
d \log (a-c) = (u - 2k) \, \omega_1 - (v - 2h) \, \omega_2.
\label{eqII7}
\end{equation}
Define the
variable $J$ to be
\be
J = \frac{1}{2} (a-c) >0.
\label{eqII8}
\end{equation}
It will appear frequently in what follows.
Equation \eqref{eqII7} then takes the form
\be
d \log J = (u - 2 k) \omega_1 - (v - 2 h) \omega_2.
\label{eqII9}
\end{equation}
The $\omega_i$ constitute a linearly independent set. Two related
coframes called $\vartheta_i$ and $\alpha_i$ can be defined in terms
of the $\omega_i$ and the functions $u$ and $v$ as follows,
\be
\begin{array}{cc}
\vartheta_1 = u \omega_1 +v \omega_2,  &  \vartheta_2 =- v \omega_1 + u \omega_2,  \\
    &       \\
\alpha_1 =u \omega_1 - v \omega_2,     &  \alpha_2 = v \omega_1 + u \omega_2.  \\
\end{array}
\label{eqII10}
\end{equation}
These relations imply that $\vartheta_1 =0$ is tangent to the level
curves specified by $H$ equals constant and $\alpha_1 =0$ is its
symmetry with respect to the principal directions.

Squaring both sides of the relation $2 H = a+c$ and subtracting
the relation $4 K = 4 ac$ yields $4 (H^2 - K) = (a-c)^2$. 
The Hodge operator, denoted by $*$,  will play an important role throughout.  
It produces the following result on the basis forms $\omega_i$,
\be
* \omega_1 = \omega_2,
\qquad
* \omega_2 =- \omega_1,
\qquad
*^2 =-1.
\label{eqII11}
\end{equation}
Moreover, adding the expressions for $da$ and $dc$ given in \eqref{eqII6},
there results
\be
\frac{1}{a - c} (da + dc) = (u-k) \omega_1 + h \omega_2 +(v-k) \omega_2 + k \omega_1
= u \omega_1 + v \omega_2 = \vartheta_1.
\label{eqII12}
\end{equation}
Finally, note that
\be
\alpha_1 + 2 \, * \omega_{12} =
u \omega_1 - v \omega_2 + 2 \, * ( h \omega_1 + k \omega_2)
= (u - 2k) \omega_1 - (v -2 h) \omega_2 = d \log J.
\label{eqII13}
\end{equation}
Therefore, the Codazzi equations \eqref{eqII12} and \eqref{eqII13}
can be summarized using the definitions of $H$ and $J$ as
\be
d H = J \vartheta_1 ,
\qquad
d \log J = \alpha_1 + 2 \, * \omega_{12}.
\label{eqII14}
\end{equation}

\section{A Theorem of Bonnet}

Suppose that $M^*$ is a surface which is isometric to $M$ such that the
principal curvatures are preserved {\bf [10-12]}. Denote all quantities which pertain
to $M^*$ with the same symbols but with asterisks, as for example
$$
a^* = a,
\qquad
c^* = c.
$$
The same notation will be applied to the variables and forms which
pertain to $M$ and $M^*$. When $M$ and $M^*$ are isometric, the
forms $\omega_i$ are related to the $\omega_i^*$ by the following transformation
\be
\omega_1^* = \cos \tau \, \omega_1 - \sin \tau \, \omega_2,
\qquad
\omega_2^* = \sin \tau \, \omega_1 + \cos \tau \, \omega_2.
\label{eqIII1}
\end{equation}

{\bf Theorem 3.1} Under the transformation of coframe given by
\eqref{eqIII1}, the associated connection forms are related by
\be
\omega_{12}^* = \omega_{12} - d \tau.
\label{eqIII2}
\end{equation}

{\bf Proof:} Exterior differentiation of $\omega_1^*$ produces
$$
d \omega_1^* = - \sin \tau \, d \tau \wedge \omega_1 + \cos \tau \, d \omega_1
- \cos \tau \, d \tau \wedge \omega_2 - \sin \tau \, d \omega_2
$$
$$
= d \tau \wedge (- \sin \tau \omega_1 - \cos \tau \omega_2 ) + \cos \tau \,
\omega_{12} \wedge \omega_2 - \sin \tau \omega_1 \wedge \omega_{12}
= (- d \tau + \omega_{12} ) \wedge \omega_2^*.
$$
Similarly, differentiating $\omega_2^*$ gives
$$
d \omega_2^* = \cos \tau \wedge \omega_1 + \sin \tau \, d \omega_1 
- \sin \tau \, d \tau \wedge \omega_2 + \cos \tau \, d \omega_2
$$
$$
= d \tau \wedge ( \cos \tau \omega_1 - \sin \tau \omega_2) +
\sin \tau \omega_{12} \wedge \omega_2 + \cos \tau \omega_1
\wedge \omega_{12} = \omega_1^* \wedge (- d \tau + \omega_{12}).
$$
$\Box$

There is a very important result which can be developed at
this point. In the case that $a= a^*$ and $c=c^*$, the 
Codazzi equations imply that
$$
\alpha_1 + 2 * \omega_{12} = d \log (a-c) = d \log ( a^* - c^*)
= \alpha_1^* + 2 * \, \omega_{12}^*.
$$
Apply the operator $*$ to both sides of this equation,
we obtain
$$
\alpha_2 - 2 \omega_{12} = \alpha_2^* - 2 \omega_{12}^*.
$$
Substituting for $\omega_{12}^*$ from Theorem 3.1, this is
\be
2 d \tau = \alpha_2 - \alpha_2^*.
\label{eqIII3}
\end{equation}

{\bf Lemma 3.1} 
$$
\vartheta_1 = \vartheta_1^*.
$$

{\bf Proof:} This can be shown in two ways. 
First from \eqref{eqIII1}, express the
$\omega_i$ in terms of the $\omega_i^*$ 
\be
\omega_1 = \cos \tau \, \omega_1^* + \sin \tau \, \omega_2^*,
\qquad
\omega_2 =- \sin \tau \, \omega_1^* + \cos \tau \, \omega_2^*.
\label{eqIII4}
\end{equation}
Therefore,
$$
\vartheta_1 = u \omega_1 + v \omega_2 =u  (\cos \tau \, \omega_1^*
+ \sin \tau \omega_2^*) + v ( - \sin \tau \omega_1^* + \cos \tau
\omega_2^*) = u^* \omega_1^* + v^* \omega_2^* = \vartheta_1^*.
$$
where $u^* = u \cos \tau - v \sin \tau$ and $v^* =u \sin \tau + 
v \cos \tau$. $\Box$

Lemma 3.1 also follows from the fact that $d H = d H^*$ and  
\eqref{eqII8}.

{\bf Lemma 3.2} 
$$
\alpha_2^* = \sin (2 \tau) \, \alpha_1 +
\cos (2 \tau) \, \alpha_2.
$$

{\bf Proof:} 
$$
\alpha_2^* = (u \sin \tau + v \cos \tau)(\cos \tau \omega_1
- \sin \tau \omega_2) + ( u \cos \tau - v \sin \tau)
( \sin \tau \omega_1 + \cos \tau \omega_2)
$$
$$
= (u \sin (2 \tau) + v \cos (2 \tau)) \omega_1 
+ (-v \sin(2 \tau) + u \cos(2 \tau)) \omega_2
$$
$$
= \sin(2 \tau) \alpha_1 + \cos(2 \tau) \alpha_2.
$$
$\Box$

Substituting $\alpha_2^*$ from Lemma 3.2 into \eqref{eqIII3},
$d \tau$ can be written as
\be
d \tau = \frac{1}{2} ( \alpha_2 - \sin(2 \tau) \alpha_1
- \cos (2 \tau) \alpha_2) = \frac{1}{2} ((1 - \cos (2 \tau))\alpha_2
- \sin (2 \tau) \alpha_1 ).
\label{eqIII5}
\end{equation}
Introduce the new variable $t = \cot (\tau)$ so $dt =- \csc^2 (\tau) \, d \tau$
and $\sin \tau = \dss\frac{1}{\sqrt{1+ t^2}}$, $\cos \tau = \dss\frac{1}{\sqrt{1+t^2}}$,
hence the following lemma.

{\bf Lemma 3.3} 
$$
dt = t \alpha_1 - \alpha_2.
$$

This is the total differential equation which must be satisfied by the
angle $\tau$ of rotation of the principal directions during the deformation.
If the deformation is to be nontrivial, it must be that this equation is
completely integrable. 

{\bf Theorem 3.2} A surface $M$ admits a non-trivial isometric deformation
that keeps the principal curvatures fixed if and only if
\be
d \alpha_1=0,  \qquad
d \alpha_2 = \alpha_1 \wedge \alpha_2
\label{eqIII6}
\end{equation}
or $\alpha_{12} = \alpha_2$.

{\bf Proof:} Differentiating both sides of Lemma 3.3 gives
$$
dt \wedge \alpha_1 + t d \alpha_1 - d \alpha_2
= ( t \alpha_1 - \alpha_2) \wedge \alpha_1 + t d \alpha_1 - d \alpha_2=0.
$$
Equating the coefficients of $t$ to zero gives the result \eqref{eqIII6}. $\Box$

This theorem seems to originate with Chern {\bf [6]} and is very
useful because it gives the exterior derivatives of the $\alpha_i$.
When the mean curvature is constant, $d H=0$, hence $d \vartheta_1 =0$,
which implies that $u=v=0$, and so $\alpha_1$ and
$\alpha_2$ vanish. Hence $dt=0$ which implies that, since the
$\alpha_i$ are linearly independent, $t$ equals a constant.
Thus, we arrive at a theorem originally due to Bonnet.

{\bf Theorem 3.3} A surface of constant mean curvature can be
isometrically deformed preserving the principal curvatures.
During the deformation, the principal directions rotate by
a fixed angle.

\section{Connection Form Associated to a Coframe and Transformation Properties}

Given the linearly independent one-forms $\omega_1, \omega_2$, the
first two of the structure equations uniquely determine the form
$\omega_{12}$. The $\omega_1, \omega_2$ are called the orthonormal
coframe of the metric
$$
ds^2 = \omega_1^2 + \omega_2^2,
$$
and $\omega_{12}$ the connection form associated with it.

{\bf Theorem 4.1} Suppose that $A>0$ is a function on $M$. 
Under the change of coframe
\be
\omega_1^* = A \omega_1,
\qquad
\omega_2^* = A \omega_2,
\label{eqIV1}
\end{equation}
the associated connection forms are related by
\be
\omega_{12}^* = \omega_{12} + \, * \, d \log A.
\label{eqIV2}
\end{equation}

{\bf Proof:} The structure equations for the transformed system
are given as
$$
d \omega_1^* = \omega_{12}^* \wedge \omega_2^*,
\qquad
d \omega_2^* = \omega_1^* \wedge \omega_{12}^*.
$$
Using \eqref{eqIV1} to replace the $\omega_i^*$ in these,
we obtain
$$
d \log \, A \wedge \omega_1 + d \omega_1 = \omega_{12}^* \wedge \omega_2,
\qquad
d \log A \wedge \omega_2 + d \omega_2 = \omega_1 \wedge \omega_{12}^*.
$$
The $\omega_i$ satisfy a similar system of structure equations,
so replacing $ d \omega_i$ here yields
$$
( \omega_{12}^* - \omega_{12}) \wedge \omega_2 = d \log \,A \wedge \omega_1,
\qquad
( \omega_{12}^* - \omega_{12}) \wedge \omega_1 =- d \log \,A \wedge \omega_2.
$$
Since the forms $\omega_i$ satisfy the equations $* \omega_1 = \omega_2$ and
$* \omega_2 =- \omega_1$, substituting these relations into the above
equations and using $\Omega_k \wedge (* \Theta_k)= \Theta_k \wedge (* \Omega_k)$,
we obtain that
in the form
$$
\omega_1 \wedge \, * ( \omega_{12}^* - \omega_{12}) =- \omega_1 \wedge d \, \log A,
\qquad
\omega_2 \wedge *( \omega_{12}^* - \omega_{12}) =- \omega_2 \wedge \, d \log \, A.
$$
Cartan's lemma can be used to conclude from these that there exist functions
$f$ and $g$ such that
$$
*( \omega_{12}^* - \omega_{12}) =- d \log \, A - f \omega_1,
\qquad
*( \omega_{12}^* - \omega_{12}) =- d \log \, A + g \omega_2.
$$
Finally, apply $*$ to both sides and use $*^2 =-1$ to obtain
$$
\omega_{12}^* - \omega_{12} =* d \log A + f \omega_2,
\qquad
\omega_{12}^* - \omega_{12} = * d \log A + g \omega_1.
$$
The forms $\omega_i$ are linearly independent, so for these two equations 
to be compatible, it suffices to put $f=g=0$, and the result follows.
$\Box$

For the necessity in the Chern criterion, Theorem 3.2, no mention of the
set $V$ of critical points of $H$ is needed. In fact, when $H$ is
constant, this criterion is met and the sufficiency also holds
with $\tau$ constant. However, when $H$ is not identically constant,
we need to take the set $V$ of critical points into account
for the sufficiency. In this case, $M-V$ is also an open, dense and
connected subset of $M$. On this subset $J>0$ and the function
$A$ can be defined in terms of the functions $u$ and $v$ as
\be
A =+ \sqrt{ u^2 + v^2} >0.
\label{eqIV3}
\end{equation}
To define more general transformations of the $\omega_i$,
define the angle $\psi$ as
\be
u = A \cos \, (\psi),
\qquad
v= A \sin (\psi).
\label{eqIV4}
\end{equation}
This angle which is defined modulo $2 \pi$, is continuous only locally and
could be discontinuous in a non-simply connected region of $M-V$. 
With $A$ and $\psi$ related to $u$ and $v$ by \eqref{eqIV4}, 
the forms $\vartheta_i$ and $\alpha_i$
can be written in terms of $A$ and $\psi$ as
\be
\begin{array}{cc}
\vartheta_1 = A ( \cos (\psi) \, \omega_1 + \sin(\psi) \omega_2),
&
\vartheta_2 = A (- \sin(\psi) \, \omega_1 + \cos(\psi) \, \omega_2),  \\
  &   \\
\alpha_1 = A (\cos(\psi) \, \omega_1 - \sin(\psi) \, \omega_2),
&
\alpha_2 = A ( \sin(\psi) \, \omega_1 + \cos(\psi) \, \omega_2).  \\
\end{array}
\label{eqIV5}
\end{equation}

The forms $\omega_i$, $\vartheta_i$, $\alpha_i$ define the same
structure on $M$ and we let $\omega_{12}$, $\vartheta_{12}$,
$\alpha_{12}$ be the connection forms associated to the 
coframes $\omega_1, \omega_2$; $\vartheta_1, \vartheta_2$;
$\alpha_1, \alpha_2$. The next Theorem is crucial for what
follows.

{\bf Theorem 4.2}
\be
\vartheta_{12} = d \psi + \omega_{12} + * d \, \log \, A = 
2 d \psi + \alpha_{12}.
\label{eqIV6}
\end{equation}

{\bf Proof:} Each of the transformations which yield the
$\vartheta_i$ and $\alpha_i$ in the form \eqref{eqIV5}
can be thought of as a composition of the two transformations which
occur in the Theorems 3.1 and 4.1. First apply the transformation
$\omega_i \raro A \omega_i$ and $\tau \raro - \psi$ with
$\omega_i^* \raro \vartheta_i$ in \eqref{eqIII1}, we get
the $\vartheta_i$ equations in \eqref{eqIV5}. Invoking Theorems 3.1 and 4.1
in turn, the first result is obtained
$$
\vartheta_{12}= d \psi + \omega_{12} + * \, d \log \, A.
$$
The transformation to the $\alpha_i$ is exactly similar except that
$\tau \raro \psi$, hence
$$
\alpha_{12} =- d \psi + \omega_{12} + \, * \, d \log \, A.
$$
This implies $* d \log \, A = \alpha_{12} + d \psi - \omega_{12}$.
When replaced in the first equation of \eqref{eqIV6}, the second 
equation appears.
Note that from Theorem 3.2, $\alpha_{12}= \alpha_2$, so the
second equation can be given as $\vartheta_{12} =2 d \psi
+ \alpha_2$.
$\Box$

Differentiating the second equation in \eqref{eqII14} and
using $d \alpha_1=0$, it follows that
\be
d \, * \omega_{12} =0.
\label{eqIV7}
\end{equation}

{\bf Lemma 4.1} The angle $\psi$ is a harmonic function
$d * d \psi =0$ and moreover, $d * \vartheta_{12} =0$.

{\bf Proof:} From Theorem 4.2, it follows by applying $*$ 
through \eqref{eqIV6} that
\be
* \vartheta_{12} = * \omega_{12} + * d \psi - d \log A
= 2 * d \psi - \alpha_1.
\label{eqIV8}
\end{equation}
Exterior differentiation of this equation using $d * \omega_{12}=0$
immediately gives
$$
d *d \psi=0.
$$
This states that $\psi$ is a harmonic function. 
Equation \eqref{eqIV8}
also implies that $d * \vartheta_{12} =0$.
$\Box$

\section{Construction of the Closed Differential Ideal Associated with M}

Exterior differentiation of the first equation in \eqref{eqII14} and
using the second equation produces
\be
d \vartheta_1 + ( \alpha_1 + 2 * \omega_{12}) \wedge \vartheta_1 =0.
\label{eqV1}
\end{equation}
The structure equation for the $\vartheta_i$ will be needed,
\be
d \vartheta_1 = \vartheta_{12} \wedge \vartheta_2 =-* \vartheta_{12} 
\wedge \vartheta_1.
\label{eqV2}
\end{equation}
From the second equation in \eqref{eqIV6}, we have $* \omega_{12} 
- d \, \log A + \alpha_1 = * d \psi$, and putting this in the first
equation of \eqref{eqIV6} we find
\be
- * \vartheta_{12} + \alpha_1 + 2 \, * \omega_{12} = 2 \, d \,\log A.
\label{eqV3}
\end{equation} 
Using \eqref{eqV3} in \eqref{eqV2}, 
\be
d \vartheta_1 + ( \alpha_1 + 2 * \omega_{12}) \wedge \vartheta_1
= 2 \, d \, \log A \wedge \vartheta_1.
\label{eqV4}
\end{equation}
Replacing $d \vartheta_1$ by means of \eqref{eqV1} implies
the following important result
\be
d \log \, A \wedge \vartheta_1 =0.
\label{eqV5}
\end{equation}
Equation \eqref{eqV5} and Cartan's lemma imply that there exists
a function $B$ such that
\be
d \, \log \, A = B \vartheta_1.
\label{eqV6}
\end{equation}
This is the first in a series of results which relates
many of the variables in question such as $J$, $B$ and $\vartheta_{12}$
directly to the one-form $\vartheta_1$. To show this requires 
considerable work. The way to proceed is to use the forms
$\alpha_i$ in Theorem 3.2 because their exterior derivatives
are known. For an arbitrary function on $M$, define
\be
df = f_1 \alpha_1 + f_2 \alpha_2.
\label{eqV7}
\end{equation}
Differentiating \eqref{eqV7} and extracting the coefficient of
$\alpha_1 \wedge \alpha_2$, we obtain
\be
f_{21}- f_{12} + f_2 =0.
\label{eqV8}
\end{equation}
In terms of the $\alpha_i$, $* d \psi = \psi_1 \alpha_2 - \psi_2 \alpha_1$,
Lemma 4.1 yields
\be
\psi_{11} + \psi_{22} +\psi_1 =0.
\label{eqV9}
\end{equation}
Finally, since $* \vartheta_{12} = 2 * d \psi - \alpha_1$, substituting
for $* d \psi$, we obtain that
\be
* \vartheta_{12} =- ( 2 \psi_2 +1) \alpha_1 + 2 \psi_1 \, \alpha_2.
\label{eqV10}
\end{equation}
Differentiating structure equation \eqref{eqV2} and using Lemma 4.1,
$$
* \vartheta_{12} \wedge \, d \vartheta_1 =0,
$$
so,
$$
* \vartheta_{12} \wedge \vartheta_{12} \wedge \vartheta_2 =0.
$$
This equation implies that either $\vartheta_{12}$ or $* \vartheta_{12}$
is a multiple by a function of the form $\vartheta_2$. Hence,
for some function $p$,
\be
\begin{array}{cc}
\vartheta_{12} =-p \vartheta_2,  &  * \vartheta_{12} =p \vartheta_1,  \\
                                 &                                \\
 \vartheta_{12} =p \vartheta_1,  &  * \vartheta_{12} =p \vartheta_2,  \\
\end{array}
\label{eqV11}
\end{equation}
Substituting the first line of \eqref{eqV11} back into
the structure equation, we have
\be
d \vartheta_1 =0.
\label{eqV12}
\end{equation}
The second line yields simply $d \vartheta_1 =p \vartheta_1 \wedge 
\vartheta_2$. Only the first case is examined now. Substituting
\eqref{eqV12} into \eqref{eqV1}, the following important constraint
is obtained
\be
( \alpha_1 + 2 * \omega_{12}) \wedge \vartheta_1 =0.
\label{eqV13}
\end{equation}

{\bf Theorem 5.1} The function $\psi$ satisfies the equation
\be
2 \psi_1 \cos ( 2 \psi) + ( 2 \psi_2 +1) \sin (2 \psi) =0.
\label{eqV14}
\end{equation}

{\bf Proof:} By substituting $* d \psi$ into \eqref{eqIV8}
we have
\be
* \vartheta_{12} = 2 * ( \psi_1 \alpha_1 + \psi_2 \alpha_2)
- \alpha_1 =- (2 \psi_2 +1) \alpha_1 + 2 \psi_1 \alpha_2.
\label{eqV15}
\end{equation}
Substituting \eqref{eqV15} into \eqref{eqIV6} and solving for
$* \omega_{12}$, we obtain that
$$
* \omega_{12} = * \vartheta_{12} - * d \psi + d \log \, A
=* \vartheta_{12} - * d \psi + B \vartheta_1 =* d \psi
- \alpha_1 + B \vartheta_1.
$$
This can be put in the equivalent form
\be
2 * \omega_{12} + \alpha_1 = 2 * d \psi - \alpha_1 + 2 B \vartheta_1.
\label{eqV16}
\end{equation}
Taking the exterior product with $\vartheta_1$ and using $d \psi_1$,
we get
$$
( \alpha_1 + 2 * \omega_{12}) \wedge \vartheta_1 =
( 2 * d \psi - \alpha_1) \wedge \vartheta_1
= ( 2 \psi_1 * \alpha_1 + 2 \psi_2 * \alpha_2 - \alpha_1) \wedge
\vartheta_1
$$
$$
= ( 2 \psi_1 \cos (2 \psi) + (2 \psi_2 +1) \sin (2 \psi)  )
\vartheta_2 \wedge \vartheta_1.
$$
Imposing the constraint \eqref{eqV13}, the coefficient of
$\vartheta_1 \wedge \vartheta_2$ can be equated to zero.
This produces the result \eqref{eqV14}.  $\Box$

As a consequence of Theorem 5.1, a new function $C$ 
can be introduced such that
\be
2 \psi_1 = C \sin ( 2 \psi),
\qquad
2 \psi_2 +1 =- C \cos (2 \psi).
\label{eqV17}
\end{equation}
Differentiation of each of these with respect to the
$\alpha_i$ basis, we get for $i=1,2$ that
$$
2 \psi_{1 i} = C_i \sin(2 \psi) + 2  \psi_i \, C \cos(2 \psi),
\qquad
2 \psi_{2i} =- C_i \cos (2 \psi) + 2 \psi_i \, C \sin (2 \psi).
$$
Substituting $f= \psi$ into \eqref{eqV8} and using the fact that
$\psi$ satisfies \eqref{eqV9} gives the pair of equations
$$
-C_1 \cos(2 \psi) - C_2 \sin (2 \psi) + 2 \psi_1 C \sin(2 \psi)
-  (2 \psi_2 +1) C \cos (2 \psi) -1=0,
$$
$$
C_1 \sin (2 \psi) - C_2 \cos (2 \psi) + 2  \psi_1 C \cos (2 \psi)
+  ( 2 \psi_2 +1) C \sin(2 \psi) =0.
$$
This linear system can be solved for $C_1$ and $C_2$ to get
\be
C_1 +C (2 \psi_2 +1) + \cos (2 \psi) =0,
\quad
C_2 -2 C \psi_1 + \sin(2 \psi) =0.
\label{eqV18}
\end{equation}
By differentiating each of the equations in \eqref{eqV18},
it is easy to verify that $C$ satisfies \eqref{eqV8},
namely $C_{12}-C_{21} -C_2=0$. Hence there exist harmonic 
functions which satisfy \eqref{eqV14}. The solution depends on two
arbitrary constants, the values of $\psi$ and $C$ at an initial
point.

{\bf Lemma 5.1} 
\be
d C = (C^2 -1) \vartheta_1,
\qquad
* \vartheta_{12}= C \vartheta_1.
\label{eqV19}
\end{equation}

{\bf Proof:} 
It is easy to express the $\vartheta_i$ in terms of 
the $\alpha_i$,
\be
\vartheta_1 = \cos (2 \psi) \alpha_1 + \sin (2 \psi) \alpha_2,
\qquad
\vartheta_2 =- \sin (2 \psi) \alpha_1 + \cos(2 \psi) \alpha_2.
\label{eqV20}
\end{equation}
Therefore, using \eqref{eqV17} and \eqref{eqV18}, it is easy to 
see that
$$
d C = C_1 \alpha_1 + C_2 \alpha_2 = (C^2 -1) (\cos (2 \psi) \alpha_1
+ \sin (2 \psi) \alpha_2) = (C^2 -1) \vartheta_1.
$$
Using \eqref{eqV17}, it follows that
$$
* \vartheta_{12} =- ( 2 \psi_2 +1) \alpha_1 + 2 \psi_1 \alpha_2 
= C \cos ( 2 \psi) \alpha_1 + C \sin (2 \psi) \alpha_2
= C ( \cos ( 2 \psi) \alpha_1 + \sin ( 2 \psi) \alpha_2) = C \vartheta_1.
$$
This implies that $\vartheta_{12} =- C \vartheta_2$. $\Box$

It is possible to obtain formulas for $B_1, B_2$.
Using \eqref{eqV20} in \eqref{eqV6}, the derivatives of
$\log A$ can be written down
\be
( \log A)_1 = B \cos (2 \psi),
\qquad
( \log A)_2 = B \sin ( 2 \psi).
\label{eqV21}
\end{equation}
Differentiating each of these in turn, we obtain for
$i=1,2$,
\be
(\log A)_{1 i} = B_i \cos (2 \psi) - 2 B \psi_i \sin (2 \psi),
\qquad
(\log A)_{2i} = B_i \sin (2 \psi) + 2 B \psi_i \cos ( 2 \psi).
\label{eqV22}
\end{equation}
Taking $f = \log A$ in \eqref{eqV8} produces a first equation for 
the $B_i$,
\be
B_1 \sin( 2 \psi) + 2 B \psi_1 \cos (2 \psi) - B_2 \cos( 2 \psi)
+ 2 B \psi_2 \sin (2 \psi) + B \sin (2 \psi) =0.
\label{eqV23}
\end{equation}
If another equation in terms of $B_1$ and $B_2$ can be found,
it can be solved simultaneously with \eqref{eqV23}.
There exists such an equation and it can be obtained
from the Gauss equation in \eqref{eqII4} which we put in the form
$$
d \omega_{12} =- a c \, \omega_1 \wedge \omega_2 =- a c \,
A^{-2} \, \alpha_1 \wedge \alpha_2.
$$

Solving \eqref{eqIV6} for $\omega_{12}$, we have
$$
\omega_{12} = d \psi + \alpha_2 +( \log A)_2 \alpha_1
- ( \log A)_1 \alpha_2.
$$
The exterior derivative of this takes the form,
$$
d \omega_{12} = [ 1 -(\log A)_{11} -( \log A)_{22} - (\log A)_1]
\alpha_1 \wedge \alpha_2.
$$
Putting this in the Gauss equation, 
$$
-(\log A)_{11} - (\log A)_{22} + \{ -(\log A)_1 +1 \} + a c A^{-2}=0.
$$
Replacing the second derivatives from \eqref{eqV22}, we have the
required second equation
\be
- B_1 \cos (2 \psi) - B_2 \sin (2 \psi) + B \{ 2 \psi_1 \sin( 2 \psi)
- ( 2 \psi_2 +1) \cos (2 \psi) \} +1 + ac A^{-2} =0.
\label{eqV24}
\end{equation}
Solving equations \eqref{eqV23} and \eqref{eqV24} together, the
following expressions for $B_1$ and $B_2$ are obtained
\be
B_1 + B (2 \psi_2 +1) - ( 1 + ac A^{-2}) \cos (2 \psi) =0,
\qquad
B_2 - 2 B \psi_1 - ( 1 + ac A^{-2}) \sin ( 2 \psi) =0.
\label{eqV25}
\end{equation}
Given these results for $B_1$ and $B_2$, it is easy to produce the 
following two Lemmas.

{\bf Lemma 5.2}
\be
d B = (BC +1 + ac A^{-2}) \vartheta_1,
\qquad
d \log J = (C + 2 B) \vartheta_1.
\label{eqV26}
\end{equation}

{\bf Proof:} Substituting \eqref{eqV25} into $d B$, we get
$$
d B = B_1 \alpha_1 + B_2 \alpha_2 = ( BC +1+ ac A^{-2})
( \cos (2 \psi) \alpha_1 + \sin(2 \psi) \alpha_2 )
= ( BC +1+ ac A^{-2}) \, \vartheta_1.
$$
Moreover,
$$
d \log J = \alpha_1 + 2 * \omega_{12} = \alpha_1
+ 2 (* \vartheta_{12} - * d \psi + d \log A)
= \alpha_1 + 2 * \vartheta_{12} - 2 * d \psi + 2 d \log A
$$
$$
= * \vartheta_{12} +2 d \log A = C \vartheta_1 +2 B \vartheta_1.
$$
$\Box$

{\bf Lemma 5.3} 
\be
d \psi = - \frac{1}{2} \sin (2 \psi) \vartheta_1 - \frac{1}{2} (C + \cos (2 \psi)) \vartheta_2.
\label{eqV27}
\end{equation}

{\bf Proof:} 
$$
2 d \psi = 2 \psi_1 \alpha_1 +2 \psi_2 \alpha_2 
= C \sin (2 \psi) \alpha_1 - (C \cos (2 \psi) +1) \alpha_2
$$
$$
= C \sin(2 \psi) ( \cos (2 \psi) \vartheta_1 - \sin(2 \psi) \vartheta_2 )
- (C \cos(2 \psi) +1)( \sin(2 \psi) \vartheta_1 + \cos(2 \psi) \vartheta_2)
$$
$$
=- \sin(2 \psi) \vartheta_1 -( C + \cos(2 \psi)) \vartheta_2.
$$
$\Box$

In the interests of completeness, it is important to verify
the following Theorem.

{\bf Theorem 5.2} The function $B$ satisfies \eqref{eqV8} provided
$\psi$ satisfies both \eqref{eqV9} and equation \eqref{eqV13}.

{\bf Proof:} Differentiating $B_1$ and $B_2$ given by \eqref{eqV25},
the left side of \eqref{eqV8} is found to be
$$
B_{21} -B_{12} +B_2 = 2 B_1 \psi_1 +B_2 ( 2 \psi_2 +1) + 2B
( \psi_{11}+ \psi_{22}+\psi_1) + A^{-2} ((ac)_1 \sin(2 \psi)
- (ac)_2 \sin(2 \psi))
$$
$$
-2 ac B A^{-2} ( \cos(2 \psi) \sin(2 \psi) - \sin(2 \psi) \cos(2 \psi))
+(1 + ac A^{-2}) ( 2 \psi_1 \cos(2 \psi) + ( 2 \psi_2 +1) \sin(2 \psi))
$$
$$
= 2 (1 +ac A^{-2}) ( 2 \psi_1 \cos (2 \psi) + ( 2 \psi_2 +1) \sin (2 \psi))
+ A^{-2} ((ac)_1 \sin(2 \psi) - (ac)_2 \cos (2 \psi)).
$$
To simplify this, \eqref{eqV9} has been substituted.
Using \eqref{eqV20} and $* d (ac) = (ac)_1 \alpha_2 - (ac)_2 \alpha_1$,
it follows that
$$
* d (ac) \wedge \vartheta_2 = (( ac)_1 \sin(2 \psi) - (ac)_2 \cos (2 \psi) ) \alpha_1 
\wedge \alpha_2.
$$
Note that the coefficient of $\alpha_1 \wedge \alpha_2$ in this appears
in the compatibility condition. To express it in another way,
begin by finding the exterior derivative of
$4 a c = (a+c)^2 - (a-c)^2$,
$$
4 d (ac) = 2 (a+ c)  (a-c) \vartheta_1 -2 (a-c)^2 ( \alpha_1 +2 * \omega_{12})
$$
Applying the Hodge operator to both sides of this, gives upon rearranging terms
$$
2 * \frac{d (ac)}{a-c} = (a+c) \vartheta_2 - (a-c) ( \alpha_2 - 2 \omega_{12}).
$$
Consequently, we can write
$$
- \frac{2}{(a-c)^2} \, * d (a c) \wedge \vartheta_2 
= ( \alpha_2 - 2 \omega_{12}) \wedge \vartheta_2
=- ( 2 \psi_1 \cos (2 \psi) + ( 2 \psi_2 +1) \sin(2 \psi)) \alpha_1 \wedge \alpha_2.
$$
Therefore, it must be that
$$
- (ac)_1 \sin(2 \psi) + (ac)_2 \cos (2 \psi) =- \frac{1}{2} (a-c)^2
( 2 \psi_1 \cos (2 \psi) + ( 2 \psi_2 +1) \sin (2 \psi)).
$$
It follows that when $f=B$, \eqref{eqV8} finally reduces to the form
$$
(1 + H^2 A^{-2}) [ 2 \psi_1 \cos( 2 \psi) +( 2 \psi_2 +1) \sin(2 \psi) ]=0.
$$
The first factor is clearly nonzero, so the second factor must vanish.
This of course is equivalent to the constraint \eqref{eqV13}.

\section{Intrinsic Characterization of M }

During the prolongation of the exterior differential
system, the additional variables $\psi$, $A$, $B$ and $C$ have been
introduced. The significance of the appearance of the function $C$,
is that the process terminates and the differentials of all these functions can
be computed without the need to introduce more functions. This means
that the exterior differential system has finally closed.

The results of the previous section, in particular the
lemmas, can be collected such that they justify the following.

{\bf Proposition 6.1} The differential system generated in terms of the
differentials of the variables $\psi$, $A$, $B$ and $C$ is closed. The variables $H$,$J,A,B,C$
remain constant along the $\vartheta_2$-curves so $\vartheta_1 =0$.
Hence, an isometry that preserves $H$ must map the $\vartheta_1$, $\vartheta_2$
curves onto the corresponding $\vartheta_1^*$, $\vartheta_2^*$ curves
of the associated surface $M^*$ which is isometric to $M$.   $\Box$

Along the $\vartheta_1$, $\vartheta_2$ curves, consider the normalized
frame,
\be
\zeta_1 = \cos( \psi) e_1 + \sin( \psi) e_2,
\qquad
\zeta_2 =- \sin(\psi) e_1 + \cos (\psi) e_2.
\label{eqVI1}
\end{equation}
The corresponding coframe and connection form are
\be
\xi_1 = \cos( \psi) \omega_1 + \sin( \psi) \omega_2,
\quad
\xi_2 =- \sin( \psi) \omega_1 + \cos( \psi) \omega_2,
\quad
\xi_{12} = d \psi + \omega_{12}.
\label{eqVI2}
\end{equation}
Then $\vartheta_1$ can be expressed as a multiple of
$\xi_1$ and $\vartheta_2, \vartheta_{12}$ in terms of $\xi_2$,
and the differential system can be summarized here:
$$
\vartheta_1 = A \xi_1,
\qquad
\vartheta_2 = A \xi_2,
\quad
\vartheta_{12} = \xi_{12} + * d \, \log A =- CA \xi_2,
$$
\be
d \log A = AB \xi_1, \qquad d B = A ( BC + 1 + a c A^{-2}) \xi_1,
\qquad d C = A ( C^2 -1) \xi_1,   
\label{eqVI3}
\end{equation}
$$
d H = AJ \xi_1,
\quad
d J = A J ( 2 B +C) \xi_1.
$$
The condition $d \vartheta_1 =0$ is equivalent to
$$
d A \wedge \xi_1 + A d \xi_1 =0.
$$
This implies that $d \xi_1 =0$ since $d A$ is proportional to $\xi_1$.
Also $d * \vartheta_{12}=0$ is equivalent to $d * \xi_{12}=0$.

Moreover, $d * \xi_{12}=0$ is equivalent
to the fact that the $\xi_1, \xi_2$-curves can be regarded as coordinate curves
parameterized by isothermal parameters. 
Therefore, along the $\xi_1, \xi_2$-curves,
orthogonal isothermal coordinates denoted $(s, t)$ can be 
introduced. The first fundamental form of $M$ then takes the form,
\be
I = \xi_1^2 + \xi_2^2 = E(s) ( ds^2 + dt^2).
\label{eqVI4}
\end{equation}
Now suppose we set $e (s) = \sqrt{E(s)}$, then
\be
\xi_1 = e(s) \, ds,
\qquad
\xi_2 = e(s) \, dt,
\qquad
\xi_{12} = \frac{e' (s)}{e^2 (s)} \xi_2 
= \frac{e' (s)}{e (s)} \, dt.
\label{eqVI5}
\end{equation}
This means such a surface is isometric to a surface of revolution.
Since $\psi$, $d * \xi_{12}=0$, \eqref{eqVI2} implies that $d * \omega_{12}=0$.
This can be stated otherwise as the principal coordinates are
isothermal and so $M$ is an isothermic surface.

Since $A,B,C,H$ and $J$ are functions of only the variable $s$,
this implies that $H$ and $J$, or $H$ and $K$, are constant along the $t$-curves with $s$
constant.  There is then the following result.

{\bf Proposition 6.2}
\be
d H \wedge d K =0,
\qquad
\xi_{12}= - (C+B)A \xi_2.
\label{eqVI6}
\end{equation}
This is equivalent to the statement $M$ is a Weingarten surface.

{\bf Proof:} The first result follows from the statement about the
coordinate system above. Since $\vartheta_{12}= \xi_{12} + * d \log A
=- CA \xi_2$ and $d A = A^2 B \xi_1$,
$$
\xi_{12} =-CA \xi_2 - * d \log A =- CA \xi_2 - * A^{-1} \, dA
= - CA \xi_2 - AB * \xi_1 =- (C+B) A \xi_2.
$$
Consequently, the geodesic curvature of each $\xi_2$-curve, $s$ constant,
is 
$$
\frac{e' (s)}{e^2 (s)} =- A (B+C),
$$
which is constant. $\Box$

To express the $\omega_i$ in terms of $ds,dt$, 
begin by writing $\omega_i$ in terms of the $\xi_i$ and then substituting
\eqref{eqVI5}, 
\be
\omega_1 = \cos ( \psi) e \, ds - \sin( \psi) e \, dt,
\qquad
\omega_2 = \sin (\psi) e \, ds + \cos( \psi) e \, dt.
\label{eqVI7}
\end{equation}
Subscripts $(s,t)$ denote differentiation and $H_s=H'$ is used interchangeably.
Beginning with $d H = H' \, ds$ and using \eqref{eqVI7}, we have
$$
dH = H_1 \omega_1 + H_2 \omega_2 = 
(H_1 \cos ( \psi) + H_2 \sin( \psi) ) \, e \, ds 
+ ( - H_1 \sin( \psi) + H_2 \cos (\psi)) \, e \, dt = H' \, ds
$$
Equating coefficients of differentials, this implies that
$$
H_1 e \cos( \psi) + H_2 e \sin (\psi) = H',
\qquad
- H_1 \sin( \psi) + H_2 \cos (\psi) =0.
$$
Solving this as a linear system  we obtain $H_1$, $H_2$,
\be
H_1 = \frac{H'}{e} \cos (\psi),
\qquad
H_2 = \frac{H'}{e} \sin( \psi).
\label{eqVI8}
\end{equation}
Noting that $u = H_1 /J$ and $v = H_2 / J$, using \eqref{eqVI2}
the forms
$\alpha_i$ can be expressed in terms of $ds, dt$
\be
\alpha_1 =\frac{H'}{J} ( \cos (2 \psi) \, ds - \sin(2 \psi) \, dt),
\qquad
\alpha_2 = \frac{H'}{J} ( \sin( 2 \psi) \, ds + \cos (2 \psi) \, dt).
\label{eqVI9}
\end{equation}
Substituting $\xi_1$ from \eqref{eqVI5} into $dH = A J \xi_1$,
$$
d H = H' ds = A J \xi_1 = A J \, e(s) \, ds.
$$
Therefore, $H' = AJ e >0$ and so $H(s)$ is an increasing function 
of $s$. Now define the function $Q (s)$ to be
\be
Q = \frac{H'}{J} = A \cdot e >0.
\label{eqVI10}
\end{equation}
Substituting \eqref{eqVI10} into \eqref{eqVI9}, the $\alpha_i$
are expressed in terms of $Q$ as well. The equations \eqref{eqIII6}
in Theorem 3.2 can easily be expressed in terms of $\psi$ and 
$Q$.

{\bf Theorem 6.1} Equations \eqref{eqIII6} are equivalent to the
following system of coupled equations in $\psi$ and $Q$:
\be
\sin(2 \psi) ( \log (Q))_s + 2 \cos(2 \psi)  \psi_s -2 \sin(2 \psi) 
\psi_t =0,
\qquad
\cos(2 \psi) ( \log (Q))_s - 2 \sin(2 \psi) \psi_s -2 \cos(2 \psi) 
\psi_t =Q.
\label{eqVI11}
\end{equation}
Moreover, equations \eqref{eqVI11} are equivalent to the
following first-order system
\be
\psi_s =- \frac{1}{2} Q \sin ( 2 \psi),
\qquad
\psi_t = \frac{1}{2} ( \log(Q))_s - \frac{1}{2} Q \cos ( 2 \psi).
\label{eqVI12}
\end{equation}
$\Box$

System \eqref{eqVI12} can be thought of as a type of Lax pair.
Moreover, \eqref{eqVI12} implies that $\psi$ is harmonic as well.
Differentiating $\psi_s$ with respect to $s$ and
$\psi_t$ with respect to $t$, it is clear that $\psi$ satisfies
Laplace's equation in the $(s,t)$ variables $\psi_{ss} + \psi_{tt}=0$.
This is another proof that $\psi$ is harmonic. 

{\bf Theorem 6.2} The function $Q(s)$ satisfies the following second-order
nonlinear differential equation
\be
Q'' (s) Q(s) - ( Q' (s) )^2 = Q^4 (s).
\label{eqVI13}
\end{equation}
There exists a first integral for this equation of the following form
\be
Q' (s)^2 = Q (s)^4 + \kappa Q(s)^2,
\quad
\kappa \in \mathbb R.
\label{eqVI14}
\end{equation}

{\bf Proof:} Equation \eqref{eqVI13} is just the compatibility
condition for the first-order system \eqref{eqVI12}.
The required derivatives are
$$
\psi_{st} =- \frac{Q}{2}  \cos ( 2 \psi)
( ( \log Q )_s - Q \cos (2 \psi)),  \quad
\psi_{ts} = \frac{1}{2} ( \log Q)_{ss} - \frac{1}{2} Q_s \cos (2 \psi) + 
Q \sin( 2 \psi) \psi_s.
$$
Equating derivatives $\psi_{st}= \psi_{ts}$,
the required \eqref{eqVI13} follows.

Differentiating both sides of \eqref{eqVI14} we get
\be
Q'' (s) = 2 Q (s)^3 + \kappa Q (s).
\label{eqVI15}
\end{equation}
Isolating $\kappa Q (s)$ from \eqref{eqVI14}
and substituting it into \eqref{eqVI15}, \eqref{eqVI13} appears.
$\Box$

It is important to note that the function $C$ which appears when the
differential ideal closes can be related to the function $Q$.

{\bf Corollary 6.1} 
\be
C = ( \frac{1}{Q} )'.
\label{eqVI16}
\end{equation}

{\bf Proof:} Using $\vartheta_i$ from \eqref{eqVI3} in Lemma 5.3,
in the $s,t$ coordinates
$$
2 d \psi = - \sin( 2 \psi) \, A e \, ds - (C + \cos (2 \psi)) \, A e \, dt
= \psi_s ds + \psi_t \, dt.
$$
Hence using \eqref{eqVI12}, this implies that $2 \psi_s =- \sin (2 \psi) \, A e
= - Q \sin (2 \psi)$, hence $Q = A e$. The second equation in \eqref{eqVI12}
for $\psi_t$ implies that $(C + \cos ( 2 \psi)) \, Ae = Q \, \cos ( 2 \psi)
- ( \log Q)'$. Replacing $Ae =Q$, this simplifies to the form \eqref{eqVI16}.

\section{Integrating the Lax Pair System}

It is clear that the first order equation in \eqref{eqVI12} for $Q (s)$
is separable and can be integrated. The integral depends on whether
$K$ is zero or nonzero:
\be
Q (s) = \frac{1}{\epsilon s + \gamma}, \quad K=0;
\qquad
\log ( \frac{2 ( K + \sqrt{K} \sqrt{Q^2 +K})}{Q}) = \epsilon \sqrt{K} s 
+ \gamma, \quad K \neq 0.
\label{eqVII1}
\end{equation}
Here $\epsilon = \pm1$ and $\gamma$ is the last constant of integration.
Taking specific choices for the constants, for example, $e^{\gamma}=2 \sqrt{K}$
when $K \neq 0$ and $ a =\sqrt{K}$, the set
of solutions \eqref{eqVII1} for $Q (s)$ can be summarized below.

\be
\begin{tabular}{cccc}
$Dom (s)$  &  $Q(s)$   &  $Dom (s)$  &   $Q (s)$  \\  \hline
       &               &             &       \\
$s >0$ & $\dss\frac{1}{s}$   &  $s<0$       &   $- \dss\frac{1}{s}$  \\
       &               &             &       \\
$0< s < \dss\frac{\pi}{a}$  &  $\dss\frac{a}{\sin (a s)}$ & $- \dss\frac{\pi}{a}< s <0$ & $-\dss \frac{a}{\sin (a s)}$ \\
       &               &             &       \\
$s >0$  & $\dss\frac{a}{\sinh (as)}$  &   $s<0$  &  $- \dss\frac{a}{\sinh (as)}$   \\  \hline
\end{tabular}
\label{eqVII2}
\end{equation}
It is presumed other choices of the constants can be geometrically
eliminated in favor of \eqref{eqVII2}.
The solutions \eqref{eqVII2} are then substituted back into linear system
\eqref{eqVI12}. The first equation in \eqref{eqVI12} implies that either
\be
\psi \equiv 0, \quad \mod \frac{\pi}{2};  \qquad
\frac{2 \psi_s}{\sin (2 \psi)} =- Q.
\label{eqVII3}
\end{equation}
Substitute $\psi \equiv 0$ into the second equation in \eqref{eqVI12}.
It implies that $(\log Q)_s = Q$ and $\psi = \pi/2$ gives $(\log Q)_s =-Q$.
In both cases $Q(s)$ is a solution which already appears in \eqref{eqVII2}.

For the second case in \eqref{eqVII3}, the equation can be put in
the form
$$
( \log | \tan ( \psi)|)_s =- Q.
$$
Integrating we have for some function $y(t)$ to be determined,
\be
\tan ( \psi) = e^{- \int Q(s) \, ds} \cdot y(t).
\label{eqVII4}
\end{equation}
Therefore $\tan (\psi)$ can be obtained by substituting for
$Q (s)$ for each of the three cases in \eqref{eqVII2}.
The upper sign holds for $s>0$ and the lower sign holds if 
$s<0$.

$(i)$ $Q (s) = \pm s^{-1}$, $- \int Q(s) ds = \log |s|^{\mp}$ and
\be
\tan (\psi) = s^{\mp} \cdot y(t).
\label{eqVII5}
\end{equation}

$(ii)$ $Q(s) = \pm \dss\frac{a}{\sin (as)}$, $- \int Q(s) ds 
= \log |\csc (as) - \cot(as)|^{\mp}$ and
\be
\tan ( \psi) = ( \tan ( \frac{as}{2}))^{\mp} \cdot y(t).
\label{eqVII6}
\end{equation}

$(iii)$ $Q (s) =\pm \dss\frac{a}{\sinh (as)}$, $- \int Q(s) ds = \mp
\arctanh ( e^{as})$, and
\be
\tan (\psi) = ( \tanh ( \frac{as}{2}))^{\mp} \cdot y(t).
\label{eqVII7}
\end{equation}
In case $(ii)$, if $s>0$ and $y(t) = \pm 1$ then $\psi = \pm \frac{1}{2} (as + \pi)$,
$\mod \pi$, and if $s<0$ and $y(t)= \pm1$, then $\psi =\pm \frac{1}{2} as$, $\mod \pi$.

It remains to integrate the second equation of the Lax pair \eqref{eqVI12}
using solutions for both $Q (s)$ and $\tan (\psi)$. The first case $(i)$
is not hard and will be shown explicitly here. The others can be
done, and more complicated cases are considered in the Appendix.

$(i)$ Consider $Q (s) = s^{-1}$ and $\tan ( \psi) = s^{-1} \cdot y(t)$.
The second equation in \eqref{eqVI12} simplifies considerably to
$y_t=-1$, therefore,
\be
y (t) =- (t + \sigma ),  \qquad
\tan (\psi) =- \frac{ ( t + \sigma)}{s}.
\label{eqVII8}
\end{equation}

For $Q(s) =- s^{-1}$ and $\tan ( \psi) = s \cdot y(t)$, the second
equation of \eqref{eqVI12} becomes $y_t =-y^2$, therefore,
\be
y (t) = \frac{1}{t + \sigma},  \qquad
\tan (\psi) = \frac{s}{t + \sigma}.
\label{eqVII9}
\end{equation}

\section{A Third Order Equation for $H$ and Fundamental Forms}

Since $\xi_{12} = (\log e(s))' dt$, using \eqref{eqVI5} $\omega_{12}$ can be written as
\be
\omega_{12} = \xi_{12} - d \psi = ( \log \, e(s))' \, dt - d \psi.
\label{eqVIII1}
\end{equation}
Using \eqref{eqII14} and \eqref{eqVI9} for $\alpha_1$, it follows that
$$
d \log (J) = Q ( \cos ( 2 \psi) \, ds - \sin (2 \psi) \, dt)
- 2 * ( \psi_t \, dt + \psi_s \, ds) + 2 * ( \log (e (s)))' \, dt.
$$
when the $\omega_i$ are put in the $s,t$ coordinates, using $* \omega_1
= \omega_2$, it can be stated that $*ds =dt$ and $*dt=-ds$. Consequently,
$d \log (J)$ simplifies to
\be
d \log (J) = ( Q \cos (2 \psi) + 2 \psi_t -2 ( \log (e (s)))') \, ds
+ (- Q \sin (2 \psi) -2 \psi_s) \, dt.
\label{eqVIII2}
\end{equation}
First order system \eqref{eqVI12} permits this to be written
using $e (s) = \sqrt{E (s)}$ as
\be
( \log (J))' + ( \log (E))' = ( \log (Q))'.
\label{eqVIII3}
\end{equation}
Hence there exists a constant $\tau$ independent of $s$ such that
$E \cdot J = \tau Q$ or
\be
E = \tau \frac{Q}{J} = \tau \frac{Q^2}{H'}
\label{eqVIII4}
\end{equation}
This result \eqref{eqVIII4} for $E$ is substituted into the Gauss
equation $- (( \log (E) )_{ss} +( \log (E))_{tt})=2 E (H^2 -J^2)$
giving
\be
( \log (E))'' = 2 ( \log (Q))'' - ( \log (H_s))''
=2 Q^2 - ( \frac{H''}{H'})'.
\label{eqVIII5}
\end{equation}
Therefore, the Gauss equation transforms into a third-order 
differential equation in the $s$ variable,
\be
( \frac{H''}{H'} )' + 2 \tau H = 2 Q^2 (1 + \tau \frac{H^2}{H'} ).
\label{eqVIII6}
\end{equation}
Thus a characterization of Bonnet surfaces is reached by means of
the solutions to these equations. This equation determines the
function $H (s)$ and after that the functions $ J (s)$ and $E (s)$.
Therefore, Bonnet surfaces have as first fundamental form the
expression
\be
I = E (s) ( d s^2 + d t^2 ), \qquad
E (s) = \tau \frac{Q^2 (s)}{H' (s)}.
\label{eqVIII7}
\end{equation}
Since $\psi$ is the angle from the principal axis $e_1$
to the $s$-curve with $t$ equals constant, the second fundamental
form is given by
\be
II = L \, ds^2 +2 M \, ds \, dt + N \, dt^2,
\label{eqVIII8}
\end{equation}
where $L,M,N$ are given by
$$
L=E (H+J \cos (2 \psi))= EH+ \tau Q \cos(2 \psi),
$$
\be
M = - EJ \sin(2 \psi)= - \tau Q \, \sin( 2 \psi),  
\label{eqVIII9}
\end{equation}
$$
N= E (H - J \cos (2 \psi)).
$$

\section{Appendix}

It is worth seeing how the second equation in \eqref{eqVI12} can be
integrated for cases $(ii)$ and $(iii)$. Only the case $s >0$
will be done with $Q(s)$ taken from \eqref{eqVII2}.

$(a)$ 
Differentiating $\tan(\psi)$ in \eqref{eqVII6}, we obtain that
$$
\psi_t = \dss\frac{\tan (\dss \frac{as}{2})}{\tan^2 ( \dss\frac{as}{2}) + y^2}
y_t (t).
$$
The following identities are required to simplify the result,
$$
\tan ( as ) = \dss\frac{2 \tan( \frac{as}{2})}
{ 1 - \tan^2 ( \frac{as}{2})},
\qquad
\cos (2 \psi) = \dss\frac{\tan^2 ( \frac{as}{2}) - y^2}
{\tan^2 ( \frac{as}{2}) +y^2}.
$$
Substituting $\psi_t$ into \eqref{eqVI12}, we obtain
$$
\dss\frac{2 \tan( \frac{as}{2})}{\tan^2 ( \frac{as}{2})+y^2} y_t
= -a \cot (as) - \frac{a}{\sin(as)} \frac{\tan^2 (\frac{as}{2}) - y^2}
{ \tan^2 ( \frac{as}{2})+y^2}.
$$
Simplifying this, we get
$$
\frac{4}{a} y_t =- \frac{1}{2} (1 - \tan^2 ( \frac{as}{2}))
- \frac{1}{2} ( \cot^2 ( \frac{as}{2}) -1) y^2 - \sec^2 ( \frac{as}{2})
+ \csc^2 ( \frac{as}{2}) y^2
$$
This simplifies to the elementary equation,
$$
y_t = \frac{a}{2} ( y^2 -1),  \qquad
y(t) =- \tanh ( \frac{at}{2} + \eta).
$$
Here $\eta$ is an integration constant. To summarize then,
$$
\tan ( \psi) = \tanh ( \frac{a t}{2} + \eta) \cdot
\tan ( \frac{as + \pi}{2}).
$$

$(b)$ Consider now $s>0$ and take $Q(s)$ from the last line of \eqref{eqVII2}.
Differentiating $\tan( \psi)$ from \eqref{eqVII7}, we get
$$
\psi_t = \frac{\coth( \dss \frac{as}{2})}
{1 + \coth^2 ( \dss\frac{as}{2}) y^2} y_t (t).
$$
In this case, the following identities are needed,
$$
\tanh (as) = \dss\frac{2 \tanh (\frac{as}{2})}{1 + \tanh^2 ( \frac{as}{2})},
\qquad
\cos (2 \psi) = \frac{1 - \coth^2 ( \frac{as}{2}) y^2}
{ 1 + \coth^2 ( \frac{as}{2}) y^2}.
$$
Therefore, \eqref{eqVI12} becomes
$$
2 \frac{\coth( \frac{as}{2})}{1 + \coth^2 ( \frac{as}{2}) y^2} y_t
= -a \coth (as) - \frac{a}{\sinh (as)} \dss\frac{\tanh^2 ( \frac{as}{2})-y^2}
{ \tanh^2 ( \frac{as}{2})+y^2}.
$$
This reduces to 
$$
- \frac{4}{a} y_t =(1 + \tanh^2 ( \frac{as}{2}) + \sech^2 ( \frac{as}{2}))
+ ( \coth^2 ( \frac{as}{2})+ 1 - \csch^2 ( \frac{as}{2})) y^2.
$$
or more simply,
$$
y_t =- \frac{a}{2} (1 + y^2),
\qquad
y(t) =- \tan ( \frac{at}{2} + \eta)
$$
To summarize then, it has been shown that,
$$
\tan( \psi) = \cot ( \frac{a t}{2} + \eta) \cdot \coth( \frac{as}{2}).
$$
These results apply to the case $s>0$ and similar results can be
found for the case $s<0$ as well.

\section{References}

\noindent
$[1]$ A. I. Bobenko, U. Eitner, Painlev\'e Equations in the Differential
Geometry of Surfaces, Lecture Notes in Mathematics, vol. 1753, Springer-Verlag,
Berlin (2000).   \\
$[2]$ G. Darboux, Th\'eorie des Surfaces, Partie 3, Paris, p. 304, (1894).  \\
$[3]$ S. S. Chern, Surface Theory with Darboux and Bianchi, Miscallanea 
Mathematica, 59-69, (1991).  \\
$[4]$ H. Hopf, Differential geometry in the large, part II, Lecture Notes in
Mathematics, vol. 1000, Springer, (1983).  \\
$[5]$ O. Bonnet, M\'emoire sur la theorie des surface applicables sur une surface
donne\'e, J. l'\'Ecole Pol., Paris, XLII Cahier, 72-92, (1867).   \\
$[6]$ S. S. Chern, Deformations of Surfaces preserving principal curvatures,
Differential Geometry and Complex Analysis, H. E. Rauch Memorial Volume,
Springer Verlag, 155-163, (1985).  \\
$[7]$ I. M. Roussos and G. E. Hernandez, On the Number of Distinct
Isometric Immersions of a Riemannian Surface into $R^3$ with given
Mean Curvature, American J. Math. {\bf 192}, 71-85, (1990).  \\
$[8]$ I. M. Roussos, Global Results on Bonnet Surfaces, J. Geometry,
{\bf 65}, 151-168, (1999).   \\
$[9]$ H. B. Lawson and R. de A. Tribuzy, On the mean curvature function for
compact surfaces, J. Diff. Geometry, {\bf 16}, 179-183, (1981).  \\
$[10]$ A. G. Colares and K. Kenmotsu, Isometric deformations of surfaces
in $R^3$ preserving the mean curvature function, Pacific J. Math.
{\bf 136, 1}, 71-80, (1989).  \\
$[11]$ X. Chen and C. K. Peng, Deformations of surfaces preserving
principal curvatures, Lecture Notes in Math., vol. 1369, Springer Verlag, 
63-73, 1987.   \\
$[12]$ S. S. Chern, W. H. Chen, K.S. Lam, Lectures in Differential
Geometry, Series on University Mathematics, vol. 1, World Scientific,
Singapore, 1999.  \\
$[13]$ K. Kenmotsu, An intrinsic characterization of $H$-deformable
surfaces, J. London Mathematical Society, {\bf 49, 2}, 555-568, (1994).  \\
$[14]$ P. Bracken, Cartan's Theory of Moving Frames and an Application
to a Theorem of Bonnet, Tensor, {\bf 70}, 261-274, (2008).  \\
\end{document}